\documentclass{amsart}
\usepackage{amsmath}
\usepackage{amssymb}
\usepackage{amsfonts}
\usepackage{amsthm}

\newtheorem{defi}{Definition}[section]
\newtheorem{teo}[defi]{Theorem}
\newtheorem{ejems}[defi]{Examples}
\newtheorem{ejem}[defi]{Example}
\newtheorem{coro}[defi]{Corollary}
\newtheorem{lem}[defi]{Lemma}
\newcommand{\prue}{\noindent \textbf{Proof:~~}}
\newcommand{\R}{\mathbb{R}}
\newcommand{\N}{\mathbb{N}} 
\newcommand{\fin}{ \hfill $\square$}
\newcommand{\tlim}{\displaystyle\lim}
\newcommand{\defrac}{\displaystyle\frac}

\begin{document}

\title[$TR$-contraction mappings on cone metric spaces]{\large  On the existence of fixed points  of contraction mappings depending of two functions on cone metric spaces}
\author[J. R. Morales and E.M. Rojas]{Jos\'e R. Morales and Edixon Rojas}
\date{}
\maketitle

\begin{center}
University of Los Andes, Faculty of Science\\
Department of Mathematics, 5101\\
 M\'erida, Venezuela \\
E-mail: moralesj@ula.ve, edixonr@ula.ve\\
\end{center}

\bigskip

\begin{abstract}
In this paper, we study the existence of fixed points for mappings defined on complete, (sequentially compact) cone metric spaces $(M,d)$, satisfying a general contractive inequality depending of two additional mappings.
\end{abstract}


\maketitle

\section{Introduction}

The concept of cone metric space was introduced by Huan Long - Guang and Zhang Xian \cite{GX}, where the set of real numbers is replaced by an ordered Banach space. They introduced the basic definitions and discuss some properties of convergence of sequences in cone metric spaces.

They also obtained various fixed point theorems for contractive single - valued maps in such spaces. Subsequently, some other mathematicians (for instance, \cite{AR,ArAzVe09,Ha09,DR,DR09,JR}) have generalized the results of Guang and Zhang \cite{GX}.

Recently, A. Beiranvand, S. Moradi, M. Omid and H. Pazandeh \cite{BMOP} introduced a new class of contractive mappings: $T-$contraction and $T-$contrative extending the Banach's contraction principle and the Edelstein's fixed point theorem, (see \cite{KK}) respectively.
Subsequently, the authors of this paper consider various extensions of classic contraction type of mappings, (more specifically: Kannan, Zamfirescu, weak-contractions and also the so-called $D(a,b)$ class)  by defining it in cone metric spaces and depending on another mapping $T$. For these classes of contractions, conditions for the existence and uniqueness of fixed points, as well for its asymptotic behavior is given \cite{MoRoA09,MoRoB09,MoRoC09}.

The purpose of this paper is to analyze the existence of fixed points for a self-map $S$ defined on a complete, (sequentially compact) cone metric space, $(M,d)$ satisfying a contraction inequality depending of two extra mappings. Our results extend some fixed points theorems of \cite{BMOP} and \cite{GX}.

\section{Preliminary facts}

Consistent with Guang and Zhang \cite{GX}, we recall the definitions of cone metric space, the notion of convergence and other results that will be needed in the sequel.

Let $E$ be a real Banach space and $P$ a subset of $E.\,\, P$ is called a \emph{cone} if and only if:
\begin{itemize}
\item[(P1)] $P$ is nonempty, closed and $P\neq \{0\};$
\item[(P2)] $a,b\in \R,\,\, a,b\geq 0$ and $x,y\in P\Rightarrow ax+by\in P;$
\item[(P3)] $x\in P$ and $-x\in P\Rightarrow x=0\Leftrightarrow P\cap (-P)=\{0\}.$
\end{itemize}
For a given cone $P\subseteq E,$ we can define a partial ordering $\leq$ on $E$ with respect to $P$ by
\begin{equation*}
x\leq y,\quad \mbox{if and only if}\quad y-x\in P.
\end{equation*}
We shall write $x<y$ to indicate that $x\leq y$ but $x\neq y,$ while $x\ll y$ will stands for $y-x\in \mbox{Int}\, P,$ where $\mbox{Int}\, P$ denotes the interior of $P.$ The cone $P\subset E$ is called \emph{normal} if there is a number $K>0$ such that for all $x,y\in E,$
$$0\leq x\leq y,\quad \mbox{implies}\quad \|x\|\leq K\|y\|.$$
The least positive number satisfying inequality above is called the \emph{normal constant} of $P$.

The cone $P$ is called \emph{regular} if every increasing sequence which is bounded from above is convergent. That is, if $(x_n)$ is a sequence such that
$$x_1\leq x_2\leq \ldots \leq x_n\leq \ldots \leq y$$
for some $y\in E,$ then there is $x\in E$ such that $\|x_n-x\|\longrightarrow 0,\,\, (n\rightarrow \infty).$

In the following we always suppose $E$ is a Banach space, $P$ is a cone with $\mbox{Int}\, P\neq \emptyset$ and $\leq$ is a partial ordering with respect to $P$.

\begin{defi}[\cite{GX}]
Let $M$ be a nonempty set. Suppose the mapping $d: M\times M\longrightarrow E$ satisfies:
\begin{itemize}
\item[(d1)] $0<d(x,y)$ for all $x,y\in M$ and $d(x,y)=0$ if and only if $x=y;$
\item[(d2)] $d(x,y)=d(y,x)$ for all $x,y\in M;$
\item[(d3)] $d(x,y)\leq d(x,z)+d(y,z)$ for all $x,y,z\in M.$
\end{itemize}
Then $d$ is called a \emph{cone metric} on $M$ and $(M,d)$ is called a \emph{cone metric space}.
\end{defi}
It is obvious that cone metric spaces generalize metric spaces.

\begin{ejems}
\begin{enumerate}
\item \textup{(\cite[Example 1]{GX})} Let $E=\R^2$, $P=\{(x,y)\in E\,:\, x,y\geq 0\}\subset \R^2,\,\, M=\R$ and $d: M\times
      M\longrightarrow E$ such that
\begin{equation*}
 d(x,y)=\bigg(|x-y|,\,\, \alpha|x-y|\bigg)
\end{equation*}
       where $\alpha\geq 0$ is a constant. Then $(M,d)$ is a cone metric space. Notice that when we consider $E$ with the Euclidean norm $\|\cdot\|_2$, the set $P$ is a normal cone, whereas if we consider $E$ with the supremum norm $\|\cdot\|_\infty$, the cone $P$ is not normal.
\item Let $E=(C[0,1],\R)$, $P=\{\varphi\in E\,:\, \varphi\geq 0\}\subset E,\,\, M=\R$ and $d:
      M\times M\longrightarrow E$ such that
\begin{equation*}
d(x,y)=|x-y|\varphi
\end{equation*}
 where $\varphi(t)=e^t\in E.$ Then $(M,d)$ is a cone metric space.
\end{enumerate}
\end{ejems}

\begin{defi}[\cite{GX}]
Let $(M,d)$ be a cone metric space. Let $(x_n)$ be a sequence in $M.$ Then:
\begin{itemize}
\item[(i)] $(x_n)$ converges to $x\in M$ if, for every $c\in E$, with $0\ll c$ there is $n_0\in \N$ such
           that for all $n\geq n_0,$

           $$d(x_n,x)\ll c.$$ We denote this by $\tlim_{n\rightarrow \infty} x_n=x$ or $x_n\longrightarrow x,\,\, (n\rightarrow \infty).$
\item[(ii)] If for any $c\in E,$ there is a number $n_0\in \N$ such that for all $m,n\geq n_0$

            $$d(x_n,x_m)\ll c,$$ then $(x_n)$ is called a \emph{Cauchy sequence} in $M;$
\item[(iii)] $(M,d)$ is a \emph{complete cone metric space} if every Cauchy sequence is convergent in $M.$
\end{itemize}
\end{defi}
The following lemma will be useful for us to prove our main results.

\begin{lem}[\cite{GX}]
Let $(M,d)$ be a cone metric space, $P$ a normal cone with normal constant $K$ and $(x_n)$ is a sequence in $M.$
\begin{itemize}
\item[(i)] $(x_n)$ converges a $x\in M$ if and only if

           $$\tlim_{n\rightarrow \infty} d(x_n,x)=0;$$
\item[(ii)] If $(x_n)$ is convergent, then it is a Cauchy sequence;

\item[(iii)] $(x_n)$ is a Cauchy sequence if and only if $\tlim_{n,m\rightarrow \infty} d(x_n,x_m)=0;$

\item[(iv)] If $x_n\longrightarrow x$ and $x_n\longrightarrow y,\,\, (n\rightarrow \infty)$ then $x=y;$

\item[(v)] If $x_n\longrightarrow x$ and $(y_n)$ is another sequence in $M$ such that
           $y_n\longrightarrow y,$ then $d(x_n,y_n)\longrightarrow d(x,y).$
\end{itemize}
\end{lem}

\begin{defi}
Let $(M,d)$ be a cone metric space. If for any sequence $(x_n)$ in $M,$ there is a subsequence $(x_{n_i})$ of $(x_n)$ such that $(x_{n_i})$ is convergent in $M$, then $M$ is called a \emph{sequentially compact cone metric space}.
\end{defi}

Next Definition and subsequent Lemma are given in \cite{BMOP} in the scope of metric spaces, here we will rewrite it in terms of cone metric spaces.
\begin{defi}

Let $(M,d)$ be a cone metric space, $P$ a normal cone with normal constant $K$ and $T: M\longrightarrow M.$ Then

\begin{itemize}
\item[(i)] $T$ is said to be \emph{continuous} if $\tlim_{n\rightarrow \infty} x_n=x$, implies that
           $\tlim_{n\rightarrow \infty} Tx_n=Tx$ for every $(x_n)$ in $M;$

\item[(ii)] $T$ is said to be \emph{sequentially convergent} if we have, for every sequence $(y_n),$ if
            $T(y_n)$ is convergent, then $(y_n)$ also is convergent;

\item[(iii)] $T$ is said to be \emph{subsequentially convergent} if we have, for every sequence $(y_n),$ if
             $T(y_n)$ is convergent, then $(y_n)$ has a convergent subsequence.
\end{itemize}
\end{defi}

\begin{lem}
Let $(M,d)$ be a sequentially compact cone metric space. Then every function $T: M\longrightarrow M$ is subsequentially convergent and every continuous function $T: M\longrightarrow M$ is sequentially convergent.
\end{lem}

\section{Main results}

In this section, first we introduce the notion of $TR-$contraction, then we extend the Banach's Contraction Principle  \cite{BMOP} and \cite{GX}.

\begin{defi}
Let $(M,d)$ be a cone metric space and $T,R,S: M\longrightarrow M$ three functions. A mapping $S$ is said to be a $TR-$\emph{contraction} if there is $a\in [0,1)$ constant such that
\begin{equation*}\label{eq3.1}
d(TSx,RSy)\leq a d(Tx,Ry)
\end{equation*} for all $x,y\in M.$
\end{defi}

\begin{ejem}\label{ej3.2}
Let $E=(C[0,1],\R)$, $P=\{\varphi\in E\,:\, \varphi\geq 0\}\subset E,\,\, M=\R$ and $d(x,y)=|x-y|e^t,$ where $e^t\in E.$ Then $(M,d)$ is a cone metric space. We consider the functions $T,R,S: M\longrightarrow M$ defined by $Tx=e^{-x}$, $Rx=2e^{-x}$ and $Sx=x+1.$ Then:
\begin{enumerate}
\item[(1)] Clearly $S$ is not a contraction;
\item[(2)] $S$ is a $TR-$contraction. In fact,
\begin{equation*}
\begin{array}{ccl}
d(TSx,RSy) &=& |TSx-RSy|e^t\\
 &=& \defrac{1}{e}|e^{-x}-2e^{-y}|e^t\\ \\
  &=& \defrac{1}{e}d(Tx,Ry)\leq ad(Tx,Ty)
\end{array}
\end{equation*}
where, $\defrac{1}{e}\leq a<1$.
\end{enumerate}
\end{ejem}

The next result extend the Theorem 1 of Guang and Zhang \cite{GX}, and Theorem 2.6 of Beiranvand, Moradi, Omid and Pazandeh \cite{BMOP}.

\begin{teo}\label{teo3.3}
Let $(M,d)$ be a complete cone metric space, $P$ be a normal cone with normal constant $K$, in addition let $T,R: M\longrightarrow M$ be  one to one and continuous functions and $S: M\longrightarrow M$ a $TR-$contraction continuous function. Then
\begin{enumerate}
\item[(i)] For every $x_0\in M$,
\begin{equation*}
\tlim_{n,m\rightarrow \infty} d(TS^nx_0,RS^{m}x_0)=0;
\end{equation*}
\item[(ii)] There exist $y_0,z_0\in M$ such that
\begin{equation}\label{eq:Thm1 01}
\tlim_{n\rightarrow \infty} TS^nx_0=y_0\quad\mbox{and}\quad \tlim_{n\rightarrow \infty} RS^nx_0=z_0;
\end{equation}
\item[(iii)] If $T$ \textup{(}or $R$\textup{)} is subsequentially convergent, then $(S^nx_0)$ has a convergent subsequence, and there exists a unique $v_0\in M$ such that
\begin{equation}
Sz_0=z_0;
\end{equation}

\item[(iv)] If $T$ \textup{(}or $R$\textup{)} is a sequentially convergent, then for each $x_0\in M$ the iterate sequence $(S^nx_0)$
      converges to $v_0$.
\end{enumerate}
\end{teo}
\prue

Let $x_0\in M$, and $(x_n)$ the Picard iteration associate to $S$  given by $x_{n+1}=Sx_n=S^nx_0$, $n=0,1,\dots$. Notice that
\begin{equation*}
d(Tx_n,Rx_{n+1})=d(TS^{n-1}x_0,RS^nx_0)\leq ad(TS^{n-2}x_0,RS^{n-1}x_0)
\end{equation*}
hence, recursively we obtain
\begin{equation}\label{eq:Thm1 02}
d(T^{n-1}x_0,RS^nx_0)\leq a^{n-1}d(Tx_0,RSx_0).
\end{equation}
Since $P$ is a normal cone with normal constant $K$, we get
\begin{equation*}
\|d(TS^{n-1}x_0,RS^nx_0)\|\leq a^{n-1}K\|d(Tx_0,RSx_0)\|
\end{equation*}
which, taking limits, implies that
\begin{equation}\label{eq:Thm1 02}
\lim_{n\to\infty}d(TS^{n-1}x_0,RS^nx_0)=0.
\end{equation}
Now, let $m,n\in\N$ with $m>n$. Then
\begin{equation*}
d(TS^{n-1}x_0,RS^{m-1}x_0)\leq d(TS^{n-1}x_0,RS^nx_0)+\cdots+ d(TS^{m-2}x_0,RS^{m-1}x_0)
\end{equation*}
if $m-n$ is odd, and
\begin{equation*}
d(TS^{n-1}x_0,RS^{m-1}x_0)\leq d(TS^{n-1}x_0,RS^nx_0)+\cdots+ d(RS^{m-2}x_0,TS^{m-1}x_0)
\end{equation*}
for $m-n$ even. Since for any $k\in\N$ can be proved analogously to \eqref{eq:Thm1 02}that
\begin{equation}\label{eq:Thm1 03}
d(RS^kx_0,TS^{k+1}x_0)\to0,\quad\mbox{as}\quad k\to\infty
\end{equation}
then, by \eqref{eq:Thm1 02} and \eqref{eq:Thm1 03}, we have that the following:
\begin{equation*}
\lim_{n,m\to\infty}d(TS^{n-1}x_0,RS^{m-1}x_0)=0
\end{equation*}
which proves (i). To prove (ii) notice
\begin{eqnarray*}
  d(TS^nx_0,TS^mx_0) &\leq& d(TS^nx_0,RS^{n+1}x_0)+d(RS^{n+1}x_0,TS^{n+2}x_0)\\
  &&+\cdots+d(RS^{m-1}x_0,TS^mx_0) \\
   &\leq&  d(TS^nx_0,RS^{n+1}x_0)+d(RS^{n+1}x_0,TS^{n+2}x_0)\\
   &&+\cdots+a^{m-2}d(RSx_0,TS^2x_0)
\end{eqnarray*}
if $m-n$ is odd, and
\begin{eqnarray*}
  d(TS^nx_0,TS^mx_0) &\leq& d(TS^nx_0,RS^{n+1}x_0)+d(RS^{n+1}x_0,TS^{n+2}x_0)\\
  &&+\cdots+d(RS^{m-1}x_0,TS^mx_0) \\
   &\leq&  d(TS^nx_0,RS^{n+1}x_0)+d(RS^{n+1}x_0,TS^{n+2}x_0)\\
   &&+\cdots+a^{m-1}d(RSx_0,TSx_0)
\end{eqnarray*}
if $m-n$ is even. As was proved above, from \eqref{eq:Thm1 02} and \eqref{eq:Thm1 03}, taking norm, considering that $P$ is a normal cone and taking limits in inequalities above, we conclude that
\begin{equation*}
\lim_{n,m\to\infty}d(TS^nx_0,TS^mx_0)=0
\end{equation*}
thus, from the fact that $(M,d)$ is a complete cone metric space,  the sequence $(TS^nx_o)$ converges. Similarly can be proved that
\begin{equation*}
\lim_{n,m\to\infty}d(RS^nx_0,RS^mx_0)=0
\end{equation*}
i.e., the sequence $(RS^nx_0)$ converges too. Therefore the limit in  \eqref{eq:Thm1 01} exist, proving in this way (ii). To prove (iii), we are going to consider that both $T$ and $R$ are subsequentially convergent. This assumption imply that $(S^nx_0)$ has a convergent subsequence. Hence, there exists $w_0\in M$ and $(n_i)_{i=1}^\infty$ such that
\begin{equation}\label{eq:Thm1 04}
\lim_{i\to\infty}S^{n_i}x_0=v_0,
\end{equation}
from the fact that $T$ and $R$ are two continuous functions, we have
\begin{equation*}
\lim_{i\to\infty}TS^{n_i}x_0=Tv_0\quad\mbox{and}\quad \lim_{i\to\infty}RS^{n_i}x_0=Rv_0
\end{equation*}
from equality \eqref{eq:Thm1 01} we conclude that
\begin{equation*}
Tv_0=y_0\quad\mbox{and}\quad Rv_0=z_0.
\end{equation*}
Since $S$ is continuous, then from \eqref{eq:Thm1 04} we get that
\begin{equation*}
\lim_{i\to\infty}S^{n_i+1}x_0=Sv_0,
\end{equation*}
also that:
\begin{equation*}
\lim_{i\to\infty}TS^{n_i+1}x_0=TSv_0\quad\mbox{and}\quad \lim_{i\to\infty}RS^{n_i+1}x_0=RSv_0.
\end{equation*}
Again by \eqref{eq:Thm1 01}, the following equalities hold
 \begin{equation*}
\lim_{i\to\infty}TS^{n_i+1}x_0=y_0\quad\mbox{and}\quad \lim_{i\to\infty}RS^{n_i+1}x_0=z_0,
\end{equation*}
hence
\begin{equation*}
TSv_0=y_0=Tv_0\quad\mbox{and}\quad RSv_0=z_0=Rv_0,
\end{equation*}
from the injectivity of $T$ and $R$ it follows that
\begin{equation*}
Sv_0=v_0.
\end{equation*}
Now, we are going to prove that the fixed point is unique. Let us suppose that another $u_0\in M$ is such that $Su_0=u_0$. Since $S$ is a $TR$-contraction, then
\begin{equation}\label{eq:Thm1 05}
d(TSv_0,RSu_0)\leq ad(Tv_0,Ru_0)
\end{equation}
but on the other hand, $d(TSv_0,RSu_0)=d(Tv_0,Ru_0)$, therefore from \eqref{eq:Thm1 05} we have that $a\geq1$ which is false. Thus the fixed point of $S$ is unique. Finally, if $T$ and $R$ are sequentially convergent, $(S^nx_0)$ is convergent and replacing $(n_i)$ by $(n)$ in  \eqref{eq:Thm1 04}, the corresponding values of the limit is $v_0$, which proves (iv).

 \fin

Taking $Tx=Rx=x$, we have the following consequence of Theorem \ref{teo3.3}:
\begin{coro}[\cite{GX}, Theorem 1]
Let $(M,d)$ be a complete cone metric space $P\subset E$ be a normal cone with normal constant $K.$ Suppose $S: M\longrightarrow M$ is a contraction function. Then $S$ has a unique fixed point in $M$ and for any $x_0\in M\,\, (S^nx)$ converges to the fixed point.
\end{coro}

Now, if we take $E=\R_{+}$ in  Theorem \ref{teo3.3} we obtain the following

\begin{coro}[\cite{BMOP}, Theorem 2.6]
Let $(M,d)$ be a complete metric space and $T: M\longrightarrow M$ be an one to one, continuous and subsequentially convergent mapping. Then  every $T-$contraction continuous function $S: M\longrightarrow M$ has a unique fixed point. Moreover, if $T$ is sequentially convergent, then for each $x_0\in M,$ the sequence $(S^nx_0)$ converge to the fixed point of $S.$
\end{coro}

If we take $E=\R$ and $Tx=Rx=x$ in  Theorem \ref{teo3.3}, then we obtain the Banach's Contraction Principle:

\begin{coro}
Let $(M,d)$ be a complete metric space and $S: M\longrightarrow M$ is a contraction mapping. Then $S$ has a unique fixed point.
\end{coro}

The following result is the localization of  Theorem \ref{teo3.3}.

\begin{teo}
Let $(M,d)$ be a complete cone metric space, $P\subset E$ be a normal cone with normal constant $K$ and $T,R: M\longrightarrow M$ be  injective, continuous and subsequentially mapping. For $c\in E$ with $0\ll c,\,\, x_0\in M,$ set
\begin{equation*}
B(Tx_0,c)=\{y\in M\,:\, d(Tx_0,y)\leq c\}.
\end{equation*}
Suppose $S: M\longrightarrow M$ is a $TR-$contraction continuous mapping for all $x,y\in B(Tx_0,c)$ and $d(TSx_0,Tx_0)\leq (1-a)c.$ Then $S$ has a unique fixed point in $B(Tx_0,c).$
\end{teo}

\prue We only need to prove that $B(Tx_0,c)$ is complete and $TSx\in B(Tx_0,c)$ for all $RSx\in B(Tx_0,c).$ Suppose that $(y_n)\subset B(Tx_0,c)$ is  a Cauchy sequence. By the completeness of $M,$ there exist $y\in M$ such that $y_n\longrightarrow y,\,\, (n\rightarrow \infty)$.

Thus, we have
\begin{equation*}
d(Tx_0,y)\leq d(y_n,Tx_0)+d(y_n,y)\leq c+d(y_n,y)
\end{equation*}
since $y_n\longrightarrow y,\,\, (n\rightarrow \infty),\,\, d(y_n,y)\longrightarrow 0.$ Hence $d(Tx_0,y)\leq c$ and $y\in B(Tx_0,c)$. Therefore, $B(Tx_0,c)$ is complete.

On the other hand, for every $Rx\in B(Tx_0,c),$
\begin{equation*}
\begin{array}{ccl}
d(Tx_0,RSx) &\leq& d(TSx_0,Tx_0)+d(TSx_0,RSx)\\
 &\leq& (1-a)c+ad(Tx_0,Rx)\leq (1-a)c+ac=c.
\end{array}
\end{equation*}
I.e., $RSx\in B(Tx_0,c)$, and so the proof is done.
\fin
\begin{coro}
Let $(M,d)$ be a complete cone metric space, $P\subset E$ be a normal cone with normal constant $K$ and $T,R: M\longrightarrow M$ be  one to one, continuous and subsequentially convergent mapping. Let suppose that $S: M\longrightarrow M$ is a mapping such that, $S^n$ is a $TR-$contraction for some $n\in \N$ and furthermore a continuous function. Then $S$ has a unique fixed point in $M.$
\end{coro}
\prue From Theorem \ref{teo3.3}, we have that $S^n$ has a unique fixed point $z_0\in M,$ that is, $S^nz_0=z_0.$ But $S^n(Sz)=S(S^nz)=Sz,$ so $S(z)$ is also fixed point of $S^n$. Hence $Sz=z$, i.e., $z$ is a fixed point of $S.$ Since the fixed point of $S$ is also fixed point of $S^n$, then the fixed point of $S$ is unique.\fin

The following example shows that we can not omit the subsequentially convergence of the function $T$ (or $R$) in Theorem \ref{teo3.3} (iii).

\begin{ejem}
Consider the Example \textup{\ref{ej3.2}}. Let $E=(C[0,1],\R),\,\, P=\{\varphi\in E\,:\, \varphi\geq 0\},\,\, M=\R$ and $d: M\times M\longrightarrow E$ defined by $d(x,y)=|x-y|e^t$ where $e^t\in E.$ Then $(M,d)$ is a complete cone metric space. Let $T,R,S: M\longrightarrow M$ be three functions defined by $Tx=e^{-x}$, $Rx=2e^{-x}$ and $Sx=x+1.$

It is clear that $S$ is a $TR-$contraction, but $T$ is not subsequentially convergent, because $Tn\rightarrow 0,\,\, (n\rightarrow \infty)$ but the sequence $(n)$ has not any convergent subsequence and $S$ has not a fixed point.\fin
\end{ejem}

In that follows by $\mathcal{F}$ we mean the family of mappings whose members are either contractive, non-expansive or $\alpha$-contraction ($0<\alpha<1$) mappings.

\begin{teo}
Let $(M,d)$ be a complete cone metric space, $P$ be a normal cone with normal constant $K$, in addition let $T,R: M\longrightarrow M$ be  one to one and continuous mappings in $\mathcal{F}$ and $S: M\longrightarrow M$ a $TR-$contraction continuous function. Then:
\begin{enumerate}
\item[(i)] For every $x_0\in M$, the iterate sequence $(S^nx_0)$ converges;
\item[(ii)] There exists a unique $v_0\in M$ such that
\begin{equation*}
Sv_0=v_0;
\end{equation*}
\item[(iii)] The iterate sequence $(S^nx_0)$ converges to the fixed point  of $S$.
\end{enumerate}
\end{teo}
\prue
(i) Let $x_0\in M$ and $(S^nx_0)$ the Picard iterate sequence
\begin{equation*}
x_{n+1}=Sx_n=S^nx_0,\quad n=0,1,\dots.
\end{equation*}
If $(S^nx_0)$ does not converges, then for each $n\in\N$
\begin{equation}\label{eq:Thm2 01}
\tlim_{n\to\infty}d(S^nx_0,S^{n+1}x_0)\nrightarrow0\qquad (n\to\infty).
\end{equation}
On the other hand, notice that
\begin{eqnarray*}
  d(TS^nx_0,RS^mx_0) &\leq& d(TS^nx_0,TS^{n+1}x_0)+d(TS^{n+1}x_0,RS^{m+1}x_0)\\
  &&+d(RS^{m+1}x_0,RS^mx_0) \\
   &\leq&d(TS^nx_0,TS^{n+1}x_0)+ad(TS^{n}x_0,RS^{m}x_0)\\
   &&+d(RS^{m+1}x_0,RS^mx_0),
\end{eqnarray*}
then
\begin{equation*}
d(TS^nx_0,RS^mx_0)\leq \frac{1}{1-a}[d(TS^nx_0,TS^{n+1}x_0)+d(RS^{m+1}x_0,RS^mx_0)]
\end{equation*}
since $T$ and $R$ are in the family $\mathcal{F}$, then inequality above can be rewrite as
\begin{equation*}
d(TS^nx_0,RS^mx_0)\leq \frac{1}{1-a}[bd(S^nx_0,S^{n+1}x_0)+cd(S^{m+1}x_0,S^mx_0)]
\end{equation*}
where $0<b,c\leq1$. By \eqref{eq:Thm2 01} we can conclude that
\begin{equation*}
\tlim_{n,m\to\infty}d(TS^nx_0,RS^mx_0)\nrightarrow0\qquad (m\to\infty),
\end{equation*}
which is a contradiction with Theorem \ref{teo3.3} (i), therefore we have that there is $v_0\in M$ such that
\begin{equation*}
\tlim_{n\to\infty}S^nx_0=v_0.
\end{equation*}
The rest of the proof runs analogous to the proof of Theorem \ref{teo3.3} with obvious changes.

\fin

\end{document}